\documentclass[a4paper,11pt]{article}
\usepackage{indentfirst,latexsym,bm,amsmath,amssymb,amsthm}
\usepackage[dvips]{graphicx}

\makeatletter\@addtoreset{equation}{section} \makeatother
 

\title{GLOBAL EXISTENCE OF THE CRITICAL SEMILINEAR WAVE EQUATIONS WITH VARIABLE
COEFFICIENTS OUTSIDE OBSTACLES}
\author{Yi Zhou
\thanks{School of Mathematical Sciences, Fudan University, Shanghai 200433, P. R. China;
Key Laboratory of Mathematics for Nonlinear Sciences (Fudan
University), Ministry of Education of China,  P. R. China; Shanghai
Key Laboratory for contemporary Applied Mathematics, School of
Mathematical Sciences, Fudan University  ({\tt
Email:yizhou@fudan.ac.cn})} \and Ning-An Lai
\thanks{School of Mathematical Sciences, Fudan University,
Shanghai 200433, P. R. China; ({\tt
Email:071018029@fudan.edu.cn})}.}
\date{}
\begin{document}
\maketitle
\begin{abstract}
In this paper, we consider exterior problem of the critical
semilinear wave
 equation in three space dimensions with variable coefficients and prove global existence of smooth solutions.
  Similar to the constant
 coefficients case, we show that the energy cannot concentrate at any
 point~$(t,x)\in(0,\infty)\times\Omega$.~For that purpose, following Ibrahim and Majdoub \cite{Ibrahim},
 we use a geometric multiplier
 close to the well-known Morawetz multiplier used in the constant coefficients
 case. Then we use
 comparison theorem from Riemannian Geometry to estimate the error
 terms. Finally, using Strichartz inequality as in Smith and Sogge \cite{Sogge}, we
 get the global existence.
 \par {\bf Keywords:}  exterior problem, variable
 coefficients wave equations, critical nonlinearity.
\end{abstract}
\section{\textbf{Introduction}}
In this paper we consider global existence of smooth solutions of
the exterior problem
\begin{equation}\label{eq13}
\left \{
\begin{aligned}
&u_{tt}-\frac{\mathrm{\partial}}{\mathrm{\partial}x_{i}}\Big(a^{ij}(x)u_{x_{j}}\Big)+u^5=0
~~~on~~(0,~\infty)\times~\Omega,\\
&u(0,x)=f(x)\in C_{0}^{\infty}(\Omega),~~~u_{t}(0,x)=g(x)\in
C_{0}^{\infty}(\Omega),\\
&u(t,x)=0~~~~~~x~\in~\partial\Omega,
\end{aligned} \right.
\end{equation}
where~$\Omega$~is the exterior of a smooth and compact obstacle
~$\vartheta\subset \mathbb{R}^3$,~$A(x)=\big(a^{ij}(x)\big)$~are
symmetric and positively definite matrices for all
~$x$~$\in~\Omega$,~$a^{ij}(x)$~are smooth functions on~$\Omega$. And
assuming the data~$(f,~g)$~satisfies a necessary compatibility
condition arising from the Dirichlet boundary condition. If
~$a^{ij}=\delta^{ij}$,~which denotes the Kronecker delta function,
we say problem (1.1) is of constant coefficients. In the case of
critical nonlinear wave equation with constant coefficients, a
wealth of results are available in the literature. For Cauchy
problem,  global existence of~$C^2$-solutions in dimension~$n=3$~was
first obtained by Rauch \cite{Rauch}, assuming the initial energy to
be small. In 1988, also for "large" data global~$C^2$-solutions in
dimension~$n=3$~were shown to exist by Struwe \cite{Struwe} in the
radially symmetric case. Grillakis \cite{Grillakis1} in 1990 was
able to remove the latter symmetry assumption and obtained the same
result. Not much later, Kapitanskii \cite{Ka} estiblished the
existence of a unique, partially regular solution for all
dimensions. Combining Strichartz inequality and Morawetz estimates,
Grillakis \cite{Grillakis2} in 1992 established global existence and
regularity for dimensions~$3\leq n\leq 5$~and announced the
corresponding results in the radial caes for dimensions~$n\leq 7$.~
Then Shatah and Struwe \cite{Shatah} obtained global existence and
regularity for dimensions~$3\leq n\leq 7$.~They also proved the
global well-posedness in the energy space in \cite{Shatah2} 1994.
For the critical exterior problem in dimension 3, Smith and Sogge
\cite{Sogge} in 1995 proved global existence of smooth solutions. In
2008, Burq et all \cite{Burq} obtained the same result in 3-D
bounded domain. \\
\indent For the critical Cauchy problem with time-independent
variable coefficients, Ibrahim and Majdoub \cite{Ibrahim} in 2003
studied the existence of both global smooth for dimensions~$3\leq n<
6$~and Shatah-Struwe's solutions for dimensions~$n\geq 3$.~\\
\indent In this paper we consider the problem (1.1) with a
general~$A(x)$~and we refer to it as critical problem with variable
coefficients. We define a
metric~$g=A^{-1}(x)=\big(a^{ij}(x)\big)^{-1},~x\in \Omega$,~then on
the Riemannian manifold~$(\Omega,~ g)$~we can introduce the distance
function~$\rho$.~To derive the global existence, the key step is to
show the~$L^6$~part of the energy associated to (1.1) cannot
concentrate at any point~$(t_{0},x_{0})$,~where~$x_{0}\in
\overline{\Omega}$.~Instead of the Morawetz
multiplier~$t\partial_{t}+r\partial_{r}+1$,~where~$r=|x|$,~we use a
geometric multiplier following Ibrahim and Majdoub [6]. That
is:~$t\partial_{t}+\rho\partial_{\rho}+1$,~ where~$\rho=\rho(x,
x_{0})$~is the distance function from some
point~$x$~to~$x_{0}$~and~$\partial_{\rho}=\nabla_{g}\rho=g^{ij}\rho_{x_{j}}\frac{\partial}
{\partial{x_{i}}}=a^{ij}\rho_{x_{j}}\frac{\partial}
{\partial{x_{i}}}$,~and~$\nabla_{g}$~here denotes the gradient on
the Riemannian manifold.
 Then we use Hessian and Laplace comparison
theorems from Riemannian Geometry to estimate the error terms.
Finally we use Strichartz estimates to obtain the global existence,
as in \cite{Sogge}.
\section{\textbf{Main result}}
In this section we show the main result and proofs.\\
Following Ibrahim and Majdoub \cite{Ibrahim}, we define:
\begin{equation}
g=A^{-1}(x)=\big(a^{ij}(x)\big)^{-1}~~~~x\in \Omega,
\end{equation}
as a Riemannian metric on~ $\Omega$,~ and consider the couple
~$(\Omega,~g)$ ~as a Riemannian manifold.
 ~For each~$x\in \Omega$,~the Riemannian metric~$g$~induces
the inner product and the norm on the tangent space
~$\Omega_{x}=\Omega$,~by:
\begin{equation}
\langle X,~Y\rangle_{g}=\langle
A^{-1}(x)X,~Y\rangle,~~~~|X|_{g}^2=\langle X,~X\rangle_{g},~~X,Y\in
\Omega,
\end{equation}
where~$\langle\cdot,\cdot\rangle$~is the standard inner product of
the Euclidean space. For~$w\in H^1(\Omega)$,~we have:
\begin{equation}
\nabla_{g}w=a^{ij}(x)w_{x_{j}} =A(x)\nabla w,~~~
|\nabla_{g}w|_{g}^2=a^{ij}(x)w_{x_{i}}w_{x_{j}}~~~x\in \Omega,
\end{equation}
where~$\nabla_{g}$~is the gradient of the Riemannian
metric~$g$,~and~$\nabla$~is the gradient on Euclidean space. Here
and in the sequence, we use geometric convention of summing
over upper and lower indices. \\
\indent In this paper we assume there are~$c_{1} > 0$~and~$c_{2} >
0$~such that
\begin{equation}
c_{1}|X|^2 \leq \big<A(x)X, ~X\big> \leq
c_{2}|X|^2,~~~~for~all~~x\in
\Omega,~~X\in \Omega,\nonumber\\
\end{equation}
then~$|\nabla_{g}w|_{g}\simeq |\nabla w|$.~\\
 \indent We define the energy of the problem (1.1):
\begin{equation}
E(t)=\frac{1}{2}\int_{\Omega}
\Big(u_{t}^2+a^{ij}(x)u_{x_{i}}u_{x_{j}}+\frac{u^6}{3}\Big)\,\mathrm{d}x.
 \end{equation}
 \subsection{\textbf{Global existence}}
\indent The key to establish global existence for (1.1) is to show
 that:
 if the data~$(f,~g)$~has compact support, and if~$u$~is a smooth
 solution to (1.1) in a half open strip
 ~$[0,~t_{0})~\times~\Omega$,~then~$u$~must be uniformly bounded
 by some constant in that strip. Then local existence and regularity
 theorems imply global existence and regularity. To establish the uniform bounds on~$u$,~by compactness it
 suffices to show that~$u$~is bounded in a neighborhood of each given point
 ~$(t_{0},x_{0})$,~where~$x_{0}\in \overline{\Omega}$.\\
{\bf Theorem 2.1.} Suppose that~$u\in C^\infty([0,~t_{0})\times
\Omega)$~ solves ~$\eqref{eq13}$.~Then if~$x_{0}\in
\overline{\Omega}$,~$u$~must be bounded in a neighborhood of
~$(t_{0},~x_{0})$,~and hence~$u\in
L^\infty([0,~t_{0})\times \Omega)$.~ \\
\indent Let us now sketch the
 proof that~$u$~cannot blow up at~$(t_{0},x_{0})$.~As in Grillakis \cite{Grillakis2} and Shatah and Struwe \cite{Shatah},
 the first key step is to show the $L^6$ part of the energy associated
 to (1.1) cannot concentrate at~$(t_{0},x_{0})$:~
\begin{equation}\label{eq311}
\lim_{t\nearrow t_{0}}\int_{\rho(x,~x_{0})\leq t_{0}-t \atop x\in
\Omega} \frac{u^6}{6}\,\mathrm{d}x=0,
 \end{equation}
 where~$\rho$~is the distance function of the matric~$g$~from~$x_{0}$~to~$x\in \Omega$.~If~$A=(\delta^{ij})$,~then~$g$~is
 the standard metric of~$\Omega$~and~$\rho(x)=|x-x_{0}|$.~For
 a general metric~$g$,~the structure of~$\rho(x)$~is more
 complicated. For the properties of this function, see section 3.\\
 \indent The proof of (2.5) will be shown in section 2.2,~following Struwe \cite{Struwe}, exploiting an
 geometric multiplier mentioned above similar to the well-known Morawetz multiplier.
  However, extra error terms
  appear in the variable case. To overcome this difficulty, we apply Hessian
  and Laplace comparison theorem from differential geometry. If~$x_{0}\in \partial \Omega$~,
  we apply the similar method as Burq et all used
 in \cite{Burq} to control the boundary term.\\
 \indent To prove Theorem 2.1,~The second key step is to use the
 Strichartz inequality to prove that $u$  is bounded near any point $(t_0,
 x_0)$,~where ~$x_{0}\in \Omega$.~ Our proof of this part is
completely parallel to Smith and Sogge \cite{Sogge}, for the
convenience of
the reader, we sketch the proofs as follows.\\
\indent Assuming identity~$\eqref{eq311}$~is hold, then combining
with the Strichartz estimates
 we show that $u\in L_{t}^4L_{x}^{12}(K)$, where $K$ is the
 domain of influence for ${(t_{0},~x_{0})}$:
\begin{equation}
K=\big\{(t,x): \rho(x,~x_{0})\leq t_{0}-t,~~(t,x)\in
[0,~t_{0})\times \Omega\big\}.
 \end{equation}
 Then Strichartz estimates shows that $u\in L_{t}^4L_{x}^{12}(K)$
 implies $\partial_{t}u\in L_{t}^\infty L_{x}^{6}(K)$. A similar argument can
 be applied to show that~$\nabla_{x}u\in L_{t}^\infty
 L_{x}^6(K)$,~which is equivalent to ~$|\nabla_{g}u|_{g}\in L_{t}^\infty
 L_{x}^6(K)$.~We then
 use H\"{o}lder's inequality to see that the total energy cannot
 concentrate at $(t_{0},~x_{0})$, that is:
 \begin{equation}
\lim_{t\nearrow t_{0}}\frac{1}{2}\int_{\rho(x,~x_{0})\leq t_{0}-t
\atop x\in \Omega}
\Big(u_{t}^2+a^{ij}(x)u_{x_{i}}u_{x_{j}}+\frac{u^6}{3}\Big)\,\mathrm{d}x=0.
 \end{equation}
\indent Now we shall give more specific details with several lemmas.
The first is to show the ~$L^6$~part of the energy cannot
concentrate at any point, the second is the spacetime estimates for
the wave equation, the third is standard and says that the energy
associated with our equation is conserved; furthermore the energy
inside spatial cross-sections of a backword
light cone is monotonic decreasing in time.\\
{\bf Lemma 2.2.} If~$u\in C^\infty([0,~t_{0})\times \Omega)$~solves
(1.1), and~$x_{0}\in \overline{\Omega}$,~then:
\begin{equation}\label{28}
\lim_{t\nearrow t_{0}}\int_{\rho(x,~x_{0})\leq t_{0}-t \atop x\in
\Omega} \frac{u^6(t,~x)}{6}\,\mathrm{d}x=0.
 \end{equation}\\
 \indent We postpone the proof of lemma 2.2 for the moment.\\
{\bf Lemma 2.3.} For the solution to the exterior problem in the
half open strip ~$[0,~t_{0})\times \Omega$:~\\
\begin{equation}\label{eq:1}
\left \{
\begin{aligned}
&u_{tt}-\frac{\mathrm{\partial}}{\mathrm{\partial}x_{i}}\Big(a^{ij}(x)u_{x_{j}}\Big)=F(t,x)
~~~on~~(0,~\infty)\times~\Omega,\\
&u(0,~x)=f(x),~~~u_{t}(0,~x)=g(x),
\end{aligned} \right.
\end{equation}
satisfies the estimates as follows:
 \begin{eqnarray}\label{210}
\|u\|_{L_{t}^{\frac{2q}{q-6}}L_{x}^q([0,~t_{0})\times~\Omega)}\leq
C\big(\|f\|_{\dot{H}^1(\Omega)}+\|g\|_{L^2(\Omega)}+
\|F\|_{L_{t}^1L_{x}^2([ 0,~t_{0})\times~\Omega)}\big)\nonumber \\
6\leq q <\infty.
 \end{eqnarray}
 \indent For the proof see Smith and Sogge \cite{Sogge}.\\
{\bf Lemma 2.4.} Let ~$u$~as above. Then
 \begin{equation}
\begin{aligned}
u\in L_{t}^{\frac{2q}{q-6}}L_{x}^q(K),~~~~if~6\leq q< \infty.\nonumber\\
\end{aligned}
 \end{equation}
{\bf Proof.} H\"{o}lder's inequality implies that if~$6<
 q<q_{1}$,~then~$L_{t}^\infty L_{x}^6\cap
 L_{t}^{\frac{2q_{1}}{q_{1}-6}}L_{x}^{q_{1}}\subset L_{t}^{\frac{2q}{q-6}}L_{x}^q$.~Since~$u\in L_{t}^\infty L_{x}^6$~
 by conservation of energy, it therefore suffices to check that
 \begin{equation}\label{2121}
\begin{aligned}
u\in L_{t}^{\frac{2q}{q-6}}L_{x}^q(K),~~~~if~10\leq q< \infty.\\
\end{aligned}
 \end{equation}
 If~$0\leq s_{1}< s_{2}< t_{0}$~, set
  \begin{equation}
\begin{aligned}
K_{s_{1}}^{s_{2}}=K\cap([s_{1}, ~s_{2}]\times \Omega),\nonumber\\
\end{aligned}
 \end{equation}
 where~$K$~is as above. Then, since~$u$~is smooth and has relatively compact
 support in~$[0, ~t_{0})\times \Omega$,~it suffices to show
 that for some fixed~$0< s_{1}< t_{0}$,~one has
  \begin{equation}
\begin{aligned}
\sup_{ s_{2}\in (s_{1},~
t_{0})}\|u\|_{L_{t}^{\frac{2q}{q-6}}L_{x}^q(K_{s_{1}}^{s_{2}})} <
\infty,~~~~
if~10\leq q< \infty.\nonumber\\
\end{aligned}
 \end{equation}
 \indent To establish this inequality, we shall want to
 apply~$\eqref{210}$,~and if the norm in the left is only taken
 over~$K_{s_{1}}^{s_{2}}$,~then the norm involving~$F$~need only be
 taken over the same set, by Huygen's principle. Thus
   \begin{equation}
\begin{aligned}
\|u\|_{L_{t}^{\frac{2q}{q-6}}L_{x}^q(K_{s_{1}}^{s_{2}})}\leq
C_{q}E_{0}+C_{q}\|u^5\|_{L_{t}^1L_{x}^2(K_{s_{1}}^{s_{2}})},\nonumber\\
\end{aligned}
 \end{equation}
 where~$E_{0}$~denotes the initial energy of ~$u$.~If~$q> 10$,~another
 application of H\"{o}lder's inequality yields
   \begin{equation}
\begin{aligned}
\|u^5\|_{L_{t}^1L_{x}^2(K_{s_{1}}^{s_{2}})}\leq \|u\|_{L_{t}^\infty
L_{x}^6(K_{s_{1}}^{s_{2}})}^{5-\frac{2q}{q-6}}\|u\|_{L_{t}^{\frac{2q}{q-6}}
L_{x}^q(K_{s_{1}}^{s_{2}})}^{\frac{2q}{q-6}},\nonumber\\
\end{aligned}
 \end{equation}
 and consequently
    \begin{equation}
\begin{aligned}
\|u\|_{L_{t}^{\frac{2q}{q-6}}L_{x}^q(K_{s_{1}}^{s_{2}})}\leq
C_{q}E_{0}+C_{q}\|u\|_{L_{t}^\infty
L_{x}^6(K_{s_{1}}^{s_{2}})}^{5-\frac{2q}{q-6}}\|u\|_{L_{t}^{\frac{2q}{q-6}}
L_{x}^q(K_{s_{1}}^{s_{2}})}^{\frac{2q}{q-6}}.\nonumber\\
\end{aligned}
 \end{equation}
 Given~$\varepsilon> 0$,~$\eqref{28}$~implies that we can choose~$s_{1}$~close enough to~$t_{0}$~so that
 \begin{equation}
 \begin{aligned}
C_{q}\|u\|_{L_{t}^\infty
L_{x}^6(K_{s_{1}}^{s_{2}})}^{5-\frac{2q}{q-6}}< \varepsilon.\nonumber\\
\end{aligned}
 \end{equation}
 If we choose
 \begin{equation}
 \begin{aligned}
\varepsilon< 2^{-\frac{2q}{q-6}}(C_{q}E_{0})^{1-\frac{2q}{q-6}},\nonumber\\
\end{aligned}
 \end{equation}
 then the following standard lemma implies that
 \begin{equation}
 \begin{aligned}
\|u\|_{L_{t}^{\frac{2q}{q-6}}L_{x}^q(K_{s_{1}}^{s_{2}})}\leq 2C_{q}E_{0},\nonumber\\
\end{aligned}
 \end{equation}
 giving us~$\eqref{2121}$~and finishing the proof.\\
{\bf Lemma 2.5.} Let~$0< C_{0}< \infty$~and suppose that~$0\leq y(s)
\in C\big([a,~ b)\big)$,~with~$y(a)=0$,~satisfies
 \begin{equation}
 \begin{aligned}
y(s)\leq C_{0}+\varepsilon y(s)^\gamma.\nonumber\\
\end{aligned}
 \end{equation}
 Then if~$\varepsilon< 2^{-\gamma}C_{0}^{1-\gamma}$~it follows that
  \begin{equation}
 \begin{aligned}
y(s)< 2C_{0},~~~~s\in [a,~ b).\nonumber\\
\end{aligned}
 \end{equation}
{\bf Proof.} Since~$C_{0}+\varepsilon x^\gamma-x<
0$~if~$\varepsilon<
 2^{-\gamma}C_{0}^{1-\gamma}$~and ~$x=2C_{0}$,~it follows that
  \begin{equation}
 \begin{aligned}
0\leq C_{0}+\varepsilon x^\gamma-x~~\forall ~x\in [0, ~x_{0}]~\Rightarrow~x_{0}< 2C_{0}.\nonumber\\
\end{aligned}
 \end{equation}
 Since~$y(s)$~must be~$\leq$~the supremum of such~$x_{0}$,~the
 lemma follows.\\
  \indent To complete the Theorem 2.1, we shall use the following
 special case of lemma 2.4:
 \begin{eqnarray*}
u\in L_{t}^4L_{x}^{12}(K).\nonumber
\end{eqnarray*}
Since~$(\partial_{t}u)_{tt}-\frac{\mathrm{\partial}}{\mathrm{\partial}x_{i}}\Big(a^{ij}(x)(\partial_{t}u)_{x_{j}}\Big)
=-5u^4\partial_{t}u$,~ and ~$\partial_{t}u$~has compact support,
estimates~$\eqref{210}\\
$~with~$q=6$~implies that, if~$0\leq s\leq
t\leq t_{0}$,~then
 \begin{equation}
\begin{aligned}
\|\partial_{t}u\|_{L_{t}^\infty L_{x}^6(K_{s}^t)}&\leq
C(s)+C\|u^4\partial_{t}u\|_{L_{t}^1L_{x}^2(K_{s}^t)}\\
&\leq
C(s)+C\|u\|_{L_{t}^4L_{x}^{12}(K_{s}^t)}^4\|\partial_{t}u\|_{L_{t}^\infty
L_{x}^6(K_{s}^t)}.\nonumber\\
\end{aligned}
 \end{equation}
 Hence if~$s$~is close enough to~$t_{0}$~so that~$C\|u\|_{L_{t}^4L_{x}^{12}}^4\leq
 \frac{1}{2}$,~we conclude that ~$\partial_{t}u\in L_{t}^\infty
L_{x}^6(K_{s}^t)$~with norm bounded by~$2C(s)$~for all~$t\in (s,~
t_{0})$,~which yields
 \begin{equation}
\begin{aligned}
\partial_{t}u\in L_{t}^\infty L_{x}^6(K).\nonumber\\
\end{aligned}
 \end{equation}
 If ~$x_{0}$~is interior to~$\Omega$,~a similar argument can be applied to show that~$\nabla_{x}u\in L_{t}^\infty
 L_{x}^6(K)$,~which implies~$\nabla_{g}u\in L_{t}^\infty
 L_{x}^6(K)$.~And from this and H\"{o}lder's inequality we
 conclude that the total energy of~$u$~cannot concentrate at~$(t_{0},~
 x_{0})$:~
  \begin{equation}\label{216}
\lim_{t\nearrow t_{0}}\frac{1}{2}\int_{\rho(x,~x_{0})\leq t_{0}-t
\atop x\in \Omega}
\Big(u_{t}^2+a^{ij}(x)u_{x_{i}}u_{x_{j}}+\frac{u^6}{3}\Big)\,\mathrm{d}x=0.
 \end{equation}
 If~$x_{0}\in \partial\Omega$,~however, this argument breaks down
 since~$\nabla_{x}u$~does not vanish on~$\partial \Omega$,~we cannot
 apply ~$\eqref{210}$~to estimate it. For the way to deal with this
 problem, one can get the details from \cite{Sogge}.\\
{\bf Lemma 2.6.} If ~$u\in C^\infty([0,~t_{0})\times \Omega)$~is a
solution to~(1.1),~then
 \begin{equation}\label{211}
\frac{1}{2}\int_{\Omega}
\Big(u_{t}^2+a^{ij}(x)u_{x_{i}}u_{x_{j}}+\frac{u^6}{3}\Big)\,\mathrm{d}x
\end{equation}
 is equal to a fixed constant ~$E_{0}<\infty$ for all $0\leq t
<t_{0}$.~ Additionally, if ~$0\leq s< t <t_{0}, x_{0}\in
\overline{\Omega}$,~ then
 \begin{eqnarray}\label{212}
\frac{1}{2}\int_{\rho(x,x_{0})\leq t_{0}-t \atop x\in \Omega}
\Big(u_{t}^2(t,x)+a^{ij}(x)u_{x_{i}}(t,x)u_{x_{j}}(t,x)+\frac{u^6(t,x)}{3}\Big)\,\mathrm{d}x\nonumber\\
\leq \frac{1}{2}\int_{\rho(x,x_{0})\leq t_{0}-s \atop x\in \Omega}
\Big(u_{t}^2(s,x)+a^{ij}(x)u_{x_{i}}(s,x)u_{x_{j}}(s,x)+\frac{u^6(s,x)}{3}\Big)\,\mathrm{d}x.
 \end{eqnarray}
{\bf Proof.} To prove the conservation of energy one multiplies both
sides of the equation
$u_{tt}-\frac{\mathrm{\partial}}{\mathrm{\partial}x_{i}}(a^{ij}(x)u_{x_{j}})+u^5=0$
by $\partial_{t}u$ to obtain the identity
 \begin{equation}\label{213}
\frac{\partial}{\partial{t}}\Big(\frac{u_{t}^2+a^{ij}(x)u_{x_{i}}u_{x_{j}}}{2}+\frac{u^6}{6}\Big)
-\frac{\mathrm{\partial}}{\mathrm{\partial}x_{i}}\Big(u_{t}a^{ij}(x)u_{x_{j}}\Big)=0.
\end{equation}\\
Thus,
 \begin{equation}
0=\frac{\partial}{\partial{t}}\int_{\Omega}\Big(\frac{u_{t}^2+a^{ij}(x)u_{x_{i}}u_{x_{j}}}{2}+\frac{u^6}{6}\Big)\,\mathrm{d}x
-\int_{\Omega}\frac{\mathrm{\partial}}{\mathrm{\partial}x_{i}}
\Big(u_{t}a^{ij}(x)u_{x_{j}}\Big)\,\mathrm{d}x.\nonumber
\end{equation}
And since the last term is always zero, by the divergence theorem,
due to the fact that ~$\partial_{t}u=0$~on~$\partial
\Omega$~and~$u(t,x)=0$~ for ~$|x|> C+t$,~we see that~$\eqref{213}$~
implies that~$\eqref{211}$~must be constant, as desired.\\
\indent To prove the other half of lemma 2.6 we need to define the
energy flux across part of the domain of dependence of a point.\\
\indent To do this, we first need to introduce some more notation.
First of all, if~$0\leq s< t <t_{0}$,~set
 \begin{equation}
K_{s}^t=K\cap ([s,~ t]\times \Omega),\nonumber
\end{equation}
where~$K$~is as above. And let~$M_{s}^t$~denote the "mantle"
associated with it:
 \begin{equation}
M_{s}^t=\partial K_{s}^t\cap ([s, ~t]\times \Omega).\nonumber
\end{equation}
Also, let~$\mathrm{d}\sigma$~denote the induced Lebesgue measure on
~$M_{s}^t$~and~$\nu=\nu(\rho,~x)=\frac{(1,~\nabla
\rho)}{\sqrt{1+|\bigtriangledown \rho|^2}}$~ denotes the unit normal
through~$(\rho,~ x)\in M_{s}^t$.~If we let~$e(u)$~be the vector
field arising from~$\eqref{213}$~,
 \begin{equation}
e(u)=\Big(\frac{u_{t}^2+a^{ij}(x)u_{x_{i}}u_{x_{j}}}{2}+\frac{u^6}{6},
-u_{t}a^{ij}(x)u_{x_{j}}\Big),\nonumber
\end{equation}
then we can define the "energy flux" across~$M_{s}^t$:~
 \begin{equation}
 \begin{aligned}
 &Flux(u,~M_{s}^t )\\
   & =\int_{M_{s}^t}\big<e(u), ~\nu \big>\mathrm{d}\sigma\\
   & =  \int_{M_{s}^t}\frac{\frac{1}{2}\big(u_{t}^2+a^{ij}(x)u_{x_{i}}u_{x_{j}}+\frac{u^6}{3}\big)-
   u_{t}a^{ij}(x)u_{x_{j}}\rho_{x_{i}}}{\sqrt{1+|\nabla \rho|^2}}
   \mathrm{d}\sigma\\
   & \geq  \int_{M_{s}^t}\frac{\frac{1}{2}\big(u_{t}^2+a^{ij}(x)u_{x_{i}}u_{x_{j}}+\frac{u^6}{3}\big)-
   \frac{1}{2}\big[u_{t}^2+\big(a^{ij}(x)u_{x_{j}}\rho_{x_{i}}\big)^2\big]}{\sqrt{1+|\nabla
   \rho|^2}}\mathrm{d}\sigma\\
   & \geq  \int_{M_{s}^t}\frac{\frac{1}{2}\big(u_{t}^2+a^{ij}(x)u_{x_{i}}u_{x_{j}}+\frac{u^6}{3}\big)-
   \frac{1}{2}\big[u_{t}^2+\big(a^{ij}(x )u_{x_{i}}u_{x_{j}}\big)\big(a^{lm}(x )\rho_{x_{l}}\rho_{x_{m}}\big)\big]}
   {\sqrt{1+|\nabla \rho|^2}}\mathrm{d}\sigma\\
   & =  \int_{M_{s}^t}\frac{\frac{1}{2}\big(u_{t}^2+a^{ij}(x)u_{x_{i}}u_{x_{j}}+\frac{u^6}{3}\big)-
   \frac{1}{2}\big(u_{t}^2+(a^{ij}(x)u_{x_{i}}u_{x_{j}})\big)}
   {\sqrt{1+|\nabla \rho|^2}}\mathrm{d}\sigma\\
   & =  \int_{M_{s}^t}\frac{u^6}{6\sqrt{1+|\nabla \rho|^2}}\mathrm{d}\sigma \geq 0,\nonumber\\
\end{aligned}
\end{equation}\\
since
~$|\nabla_{g}\rho|_{g}^2=g^{ij}(x)\rho_{x_{i}}\rho_{x_{j}}=a^{ij}(x)\rho_{x_{i}}\rho_{x_{j}}=1$.~
Also, Cauchy-Schwarz inequality is used to prove the above
inequality. If we integrate~$\eqref{213}$~over~$K_{s}^t$~we arrive
at the "flux identity":
 \begin{equation}\label{214}
 \begin{aligned}
&\frac{1}{2}\int_{\rho(x,~x_{0})\leq t_{0}-t \atop x\in \Omega}
\Big(u_{t}^2(t,~x)+a^{ij}(x)u_{x_{i}}(t,~x)u_{x_{j}}(t,~x)+\frac{u^6(t,~x)}{3}\Big)
\,\mathrm{d}x+Flux\big(u, ~M_{s}^t\big)\\
&= \frac{1}{2}\int_{\rho(x,~x_{0})\leq t_{0}-s \atop x\in \Omega}
\Big(u_{t}^2(s,~x)+a^{ij}(x)u_{x_{i}}(s,~x)u_{x_{j}}(s,~x)+\frac{u^6(s,~x)}{3}\Big)\,\mathrm{d}x,
\end{aligned}
\end{equation}
 that is
\begin{equation}
\tag{\ref{214}$'$} \label{2144}
 \begin{aligned}
E\big(u,~D(t)\big)+Flux\big(u, ~M_{s}^t\big)=E\big(u,~D(s)\big),
\end{aligned}
\end{equation}
 where
\begin{eqnarray*}
E\big(u,~D(t)\big)=\frac{1}{2}\int_{\rho(x,~x_{0})\leq t_{0}-t \atop
x\in \Omega}
\big(u_{t}^2+a^{ij}(x)u_{x_{i}}u_{x_{j}}+\frac{u^6}{3}\big)\,\mathrm{d}x.\nonumber
\end{eqnarray*}
Since~$Flux(u,~M_{s}^t)\geq 0$,~we
see~$\eqref{214}$~implies~$\eqref{212}$,~which completes the
proof.\\
\indent And we conclude from~$\eqref{2144}$~that~$t \rightarrow
E\big(u,~D(t)\big)$~is a non-increasing function on~$[0,~t_{0})$.~
It is also bounded, since~$E\big(u,~D(t)\big) \leq E(t) \leq
 E_{0} < \infty$,~on account of our assumptions on the data. Hence,~$E\big(u,~D(t)\big)$~and
 ~$E\big(u,~D(s)\big)$~in~$\eqref{2144}$~must approach a common
limit. This in turn gives the important
 fact that
   \begin{equation}\label{21777}
 \begin{aligned}
Flux(u,~M_{s}^t ) \rightarrow 0,~~~as~s \rightarrow
t.\\
\end{aligned}
\end{equation}
 \indent Given~$\varepsilon> 0$,~ from the identity~$\eqref{216}$,~we
 can find a~$0< t_{1}< t_{0}$~so that
 \begin{equation}
\frac{1}{2}\int_{\rho(x,~x_{0})\leq t_{0}-t_{1} \atop x\in \Omega}
\Big(u_{t}^2+a^{ij}(x)u_{x_{i}}u_{x_{j}}+\frac{u^6}{3}\Big)(t_{1},~x)\,\mathrm{d}x<
\frac{\varepsilon}{2}.\nonumber
 \end{equation}
 By dominated convergence, there is a~~$\delta >0$~so that
  \begin{equation}
\frac{1}{2}\int_{\rho(x,~x_{0})\leq \delta+t_{0}-t_{1} \atop x\in
\Omega}
\Big(u_{t}^2+a^{ij}(x)u_{x_{i}}u_{x_{j}}+\frac{u^6}{3}\Big)(t_{1},~x)\,\mathrm{d}x<
\varepsilon.\nonumber
 \end{equation}
 Then by the monotonicity of energy~$\eqref{212}$,~yields
 \begin{equation}
\int_{\rho(x,~x_{0})\leq \delta+t_{0}-t \atop x\in \Omega}
\frac{u^6(t,~x)}{6}\,\mathrm{d}x<
\varepsilon,~~~~~t_{1}\leq t\leq t_{0}.\nonumber\\
 \end{equation}
 \indent Let
  \begin{equation}
K^\delta=\big\{(t,~x): \rho(x,x_{0})<\delta+ t_{0}-t,~~(t,~x)\in
[0,~t_{0})\times \Omega\big\}.\nonumber\\
 \end{equation}
For~$\varepsilon$~sufficiently small, we can repeat the proof of
lemma 2.4 with~$K$~replaced by~$K^\delta$,~to conclude that
 \begin{equation}
u\in L_{t}^4L_{x}^{12}(K^\delta).\nonumber\\
 \end{equation}
 Combing with lemma 2.6 as above we can now argue as before to conclude that
  \begin{equation}
\partial_{t}u\in L_{t}^\infty L_{x}^6(K^\delta),~~\nabla_{x}u\in L_{t}^\infty L_{x}^6(K^\delta),\nonumber\\
 \end{equation}
 which implies~$u\in L^\infty(K^{\frac{\delta}{2}})$~by Sobolev's theorem. Since~$u$~vanishes outside
 of a relatively compact subset of~$[0,~t_{0})\times
\Omega$,~we can cover its support by finitely many of these
sets~$K^{\frac{\delta}{2}}$.~ Hence,~$u\in L^\infty
([0,~t_{0})\times
\Omega)$,~which implies that~$u$~can be extended to a global solution.\\
\indent For~$x_{0}\in \partial\Omega$,~an additional argument is
needed since~$\nabla_{x}u$~does not vanish on~$\partial\Omega$.~Here
we skip this step as the method is just totally the same as Smith
and Sogge used in \cite{Sogge}.
\subsection{\textbf{Nonconcentration of ~$L^6$~part of energy}}
 \indent Now we prove lemma 2.2. For that purpose we need several lemmas about
 differential geometry. And we work on~$\Omega$~with metric~$g=\langle\cdot,
 \cdot\rangle_{g}$~given by (2.1).\\
{\bf Lemma 2.7.}\label{L4} Let~$f$~be function and~$X\in
\Omega_{x}$~be vector field. Then, we have
\begin{equation}\label{217}
\big<\nabla_{g}f,~
\nabla_{g}\big(X(f)\big)\big>_{g}=\big<\nabla_{\nabla_{g}f}X,~
\nabla_{g}f\big>_{g}+X\big(\frac{1}{2}|\nabla_{g}f|_{g}^2\big),~~~x\in
\Omega.
 \end{equation}
 \indent We shall prove this identity in section 3.\\
\indent Since for any vector field~$Y, Z\in \Omega_{x}$,~we have
\begin{eqnarray*}
 \big<\nabla_{Y}\nabla_{g}(\rho^2), ~Z\big>_{g}
   & = & Y\big<\nabla_{g}(\rho^2), ~Z\big>_{g}-\big<\nabla_{g}(\rho^2), ~\nabla_{Y}Z\big>_{g}\\
   & = & YZ(\rho^2)-(\nabla_{Y}Z)(\rho^2)\\
   & = & D^2\rho^2(Y, ~Z),
 \end{eqnarray*}
 where~$D^2\rho^2$~is the Hessian of the function~$\rho^2$, if we
 replace~$X, ~f$~in the equality~$\eqref{217}$
 ~with~$\nabla_{g}(\frac{1}{2}\rho^2),~
 u$~respectively, we get
 \begin{equation}
\tag{\ref{217}$'$} \label{2177}
\begin{aligned}
 \big<\nabla_{g}u,~
\nabla_{g}\big(X(u)\big)\big>_{g} &= \big<\nabla_{g}u,~
\nabla_{g}\Big(\nabla_{g}(\frac{1}{2}\rho^2)(u)\Big)\big>_{g}
= \big<\nabla_{g}u,~ \nabla_{g}\big(\rho a^{ij}\rho_{x_{j}}u_{x_{i}}\big)\big>_{g}\\
    &=\big<\nabla_{\nabla_{g}u}\nabla_{g}\big(\frac{1}{2}\rho^2\big),~
\nabla_{g}u\big>_{g}+\nabla_{g}\big(\frac{1}{2}\rho^2\big)\big(\frac{1}{2}|\nabla_{g}u|_{g}^2\big)\\
&=\frac{1}{2}D^2\rho^2\big(\nabla_{g}u,~
\nabla_{g}u\big)+\nabla_{g}\big(\frac{1}{2}\rho^2\big)\big(\frac{1}{2}|\nabla_{g}u|_{g}^2\big).
   \end{aligned}
 \end{equation}
{\bf Lemma 2.8.} If the sectional curvature~$\kappa$~of the
Riemannian manifold~$(\Omega, ~g)$~satisfies
\begin{eqnarray*}
-a^2\leq \kappa \leq a^2
 \end{eqnarray*}
then for the distance function~$\rho$~on~$(\Omega, ~g)$,~$\forall
~X,~Y \in \Omega_{x}$,~we have
 \begin{equation}\label{eq217}
\begin{aligned}
\lim_{\rho \rightarrow 0}\frac{\partial}{\partial_{x_{i}}}\big(\rho
a^{ij}\rho_{x_{j}}\big)=\lim_{\rho \rightarrow
0}\Big(\triangle_{g}\big(\frac{1}{2}\rho^2\big)-\frac{1}{2}g^{lm}\frac{\partial
g_{lm}}{\partial x_{i}}\rho\nabla_{g}\rho\Big)=3,\\
\end{aligned}
 \end{equation}
  \begin{equation}\label{eq218}
\begin{aligned}
\lim_{\rho \rightarrow 0}\frac{1}{2}D^2\rho^2(X,~ Y)=\big<X,~ Y\big>_{g},\\
   \end{aligned}
 \end{equation}
 where~$\triangle_{g}$~is the Laplace operator on~$(\Omega,~ g)$.~ \\
 \indent For the proof see section
 3.\\
{\bf Lemma 2.9.} Assume that~$u$~is a weak solution to (1.1), then
we have
 \begin{equation}\label{eq2177}
\begin{aligned}
\Big|\Big| \frac{\partial u}{\partial \nu}\Big|\Big|_{L^2((0,
~t_{0}) \times
\partial\Omega)}
\leq CE(u)^{\frac{1}{2}},
\end{aligned}
 \end{equation}
 where~$\frac{\partial u}{\partial \nu}$~is the trace to the
 boundary of the exterior normal derivative of~$u$.~\\
{\bf Proof.} Similar to the constant case in Burq et all [1],
take~$Z \in C^{\infty}(\Omega;~T\Omega)$~a vector field whose
restriction to~$\partial\Omega$~is equal
to~$\frac{\partial}{\partial\nu}$~and compute for~$0 < T <t_{0}$~
 \begin{equation}
\begin{aligned}
&\int_{0}^{T}\int_{\Omega}\Big[\big(\partial_{t}^{2}-\frac{\partial}{\partial
x_{i}}(a^{ij}\frac{\partial}{\partial
x_{j}})\big),Z\Big]u(t,~x)\cdot u(t,~x)\,\mathrm{d}x\,\mathrm{d}t\\
=&\int_{0}^{T}\int_{\Omega}\Big[\big(\partial_{t}^{2}-a^{ij}\frac{\partial^{2}}{\partial
x_{i}\partial x_{j}}-\frac{\partial a^{ij}}{\partial
x{i}}\frac{\partial}{\partial
x{j}}\big)Zu-Z\big(\partial_{t}^{2}-a^{ij}\frac{\partial^{2}}{\partial
x_{i}\partial x_{j}}-\frac{\partial a^{ij}}{\partial
x{i}}\frac{\partial}{\partial
x{j}}\big)u\Big]u\,\mathrm{d}x\,\mathrm{d}t\nonumber.\\
\end{aligned}
 \end{equation}
 Integrate by parts, we obtain
  \begin{equation}\label{eq2178}
\begin{aligned}
&\int_{0}^{T}\int_{\Omega}\Big[\big(\partial_{t}^{2}-\frac{\partial}{\partial
x_{i}}(a^{ij}\frac{\partial}{\partial
x_{j}})\big),Z\Big]u(t,~x)\cdot u(t,~x)\,\mathrm{d}x\,\mathrm{d}t\\
=&\int_{0}^{T}\int_{\Omega}\frac{\partial}{\partial
x_{i}}\big[(Zu)a^{ij}\frac{\partial u}{\partial
x_{j}}\big]\,\mathrm{d}x\,\mathrm{d}t+\int_{0}^{T}\int_{\Omega}-(Zu)u^{5}+Z(u^{5})u\,\mathrm{d}x\,\mathrm{d}t\\
+&\Big[\int_{\Omega}\partial_{t}(Zu)\cdot
u\,\mathrm{d}x\Big]_{0}^{T}-\Big[\int_{\Omega}(Zu)\cdot
\partial_{t}u\,\mathrm{d}x\Big]_{0}^{T}.\\
\end{aligned}
 \end{equation}
 From the assumption of the coefficients~$a^{ij}$,~and noting that on~$[0,~T] \times \partial\Omega$,~$\nabla
 _{x}u=(\partial_{\nu}u)\nu$,~ we have
 \begin{equation}\label{eq2179}
\begin{aligned}
&\int_{0}^{T}\int_{\Omega}\frac{\partial}{\partial
x_{i}}\big[(Zu)a^{ij}\frac{\partial u}{\partial
x_{j}}\big]\,\mathrm{d}x\,\mathrm{d}t=\int_{0}^{T}\int_{\partial\Omega}\frac{\partial
u}{\partial \nu}a^{ij}\frac{\partial u}{\partial x_{j}}\nu_{i}\,\mathrm{d}\sigma\,\mathrm{d}t\\
&=\int_{0}^{T}\int_{\partial\Omega}\frac{\partial u}{\partial
\nu}a^{ij}(\partial_{\nu}u)\nu_{i}\nu_{j}\,\mathrm{d}\sigma\,\mathrm{d}t
\geq C\int_{0}^{T}\int_{\partial\Omega}\big|\frac{\partial
u}{\partial
\nu}\big|^{2}\,\mathrm{d}\sigma\,\mathrm{d}t.\\
\end{aligned}
 \end{equation}
 Remark now that if~$Z=\sum_{j}b_{j}\frac{\partial}{\partial
 x_{j}}$,~then integration by parts yields(using the Dirichlet boundary condition)
  \begin{equation}\label{eq21710}
\begin{aligned}
&\Big|\int_{\Omega}-(Zu)u^{5}+Z(u^{5})u\,\mathrm{d}x\,\mathrm{d}t\Big|=\frac{4}{6}\Big|\int_{0}^{T}\int_{\Omega}
Z(u^{6})(t,~x)\,\mathrm{d}x\,\mathrm{d}t\Big| \\
=&\frac{4}{6}\Big|\int_{0}^{T}\int_{\Omega} \sum_{j}\frac{\partial
b_{j}}{\partial x_{j}}u^{6}\,\mathrm{d}x\,\mathrm{d}t\Big|\leq
CE(u).\\
\end{aligned}
 \end{equation}
 while
  \begin{equation}\label{eq21711}
\begin{aligned}
\Big|\Big[\int_{\Omega}\partial_{t}(Zu)\cdot
u\,\mathrm{d}x\Big]_{0}^{T}-\Big[\int_{\Omega}(Zu)\cdot
\partial_{t}u\,\mathrm{d}x\Big]_{0}^{T}\Big|\leq CE(u).\\
\end{aligned}
 \end{equation}
 and~$\Big[\big(\partial_{t}^{2}-\frac{\partial}{\partial
x_{i}}(a^{ij}\frac{\partial}{\partial
x_{j}})\big),Z\Big]=-\big[\frac{\partial}{\partial
x_{i}}(a^{ij}\frac{\partial}{\partial x_{j}}),Z\big]$~as a second
order differential operator in the ~$x$~variable is continuous
from~$H_{0}^{1}(\Omega)$~to~$H^{-1}(\Omega)$~and consequently
  \begin{equation}\label{eq21712}
\begin{aligned}
\Big|\int_{0}^{T}\int_{\Omega}\Big[\big(\partial_{t}^{2}-\frac{\partial}{\partial
x_{i}}(a^{ij}\frac{\partial}{\partial
x_{j}})\big),Z\Big]u(t,~x)\cdot
u(t,~x)\,\mathrm{d}x\,\mathrm{d}t\Big|\leq CE(u).\\
\end{aligned}
 \end{equation}
As the constants are uniform with respect
to~$0<T<t_{0}$,~collecting~\eqref{eq2178}, \eqref{eq2179},
\eqref{eq21710}, \eqref{eq21711} and
\eqref{eq21712}~yields~$\eqref{eq2177}$.~\\
{\bf Proof of lemma 2.2.} Following Ibrahim and Majdoub
\cite{Ibrahim}, we use
 a geometric multiplier. For the sake of notation it
 is convenient to shift~$(t_{0}, ~x_{0})\in \mathbb{R} \times
 \Omega$~to the origin.\\
 \indent Multiply the
 equation~$u_{tt}-\frac{\mathrm{\partial}}{\mathrm{\partial}x_{i}}(a^{ij}(x)u_{x_{j}})+u^5=0$~
 by~$tu_{t}+\rho a^{lm}\rho_{x_{m}}u_{x_{l}}+u$.~By~$\eqref{213}$~it is easy
 to see the contribution from the first term is
 \begin{equation}
\frac{\partial}{\partial
t}\big[t(\frac{1}{2}\big(u_{t}^2+a^{ij}(x)u_{x_{i}}u_{x_{j}}+\frac{u^6}{3}\big)\big]
-\frac{1}{2}\big(u_{t}^2+a^{ij}(x)u_{x_{i}}u_{x_{j}}+\frac{u^6}{3}\big)
-\frac{\mathrm{\partial}}{\mathrm{\partial}x_{i}}\big(tu_{t}a^{ij}(x)u_{x_{j}}\big)=0.\nonumber
 \end{equation}
 Similarly, we compute
 \begin{equation}
0=\Big(u_{tt}-\frac{\mathrm{\partial}}{\mathrm{\partial}x_{i}}\big(a^{ij}(x)u_{x_{j}}\big)+u^5\Big)\Big(\rho
a^{lm}\rho_{x_{m}}u_{x_{l}}\Big).\nonumber
 \end{equation}
 For
\begin{eqnarray*}
 \rho a^{lm}\rho_{x_{m}}u_{x_{l}}u_{tt}
   & = & (\rho a^{lm}\rho_{x_{m}}u_{x_{l}}u_{t})_{t}-\rho a^{lm}\rho_{x_{m}}u_{x_{l}t}u_{t}\\
   & = & (\rho a^{lm}\rho_{x_{m}}u_{x_{l}}u_{t})_{t}-\frac{1}{2}\rho a^{lm}\rho_{x_{m}}\frac
   {\partial u_{t}^2}{\partial x_{l}}\\
   & = & (\rho a^{lm}\rho_{x_{m}}u_{x_{l}}u_{t})_{t}-\frac{1}{2}\big[\frac{\partial}{\partial x_{l}}
   (\rho a^{lm}\rho_{x_{m}}u_{t}^2)-\frac{\partial}{\partial x_{l}}(\rho
   a^{lm}\rho_{x_{m}})u_{t}^2\big].\nonumber
 \end{eqnarray*}
 Using~\eqref{2177}~in lemma 2.7, we have
\begin{equation}
\begin{aligned}
 &\rho a^{lm}\rho_{x_{m}}u_{x_{l}}\frac{\partial}{\partial
 x_{i}}\big(a^{ij}(x)u_{x_{j}}\big)\\
&=\frac{\partial}{\partial x_{i}}\big(\rho
a^{lm}\rho_{x_{m}}u_{x_{l}}a^{ij}(x)u_{x_{j}}\big)-
   \frac{\partial}{\partial x_{i}}\big(\rho a^{lm}\rho_{x_{m}}u_{x_{l}}\big)a^{ij}(x)u_{x_{j}}\\
   &=\frac{\partial}{\partial x_{i}}\big(\rho
   a^{lm}\rho_{x_{m}}u_{x_{l}}a^{ij}(x)u_{x_{j}}\big)-\big<\nabla_{g}u,~ \nabla_{g}(\rho
   a^{lm}\rho_{x_{m}}u_{x_{l}})\big>_{g}\\
   &=\frac{\partial}{\partial x_{i}}(\rho
   a^{lm}\rho_{x_{m}}u_{x_{l}}a^{ij}u_{x_{j}})-\frac{1}{2}D^2\rho^2\big(\nabla_{g}u,~
\nabla_{g}u\big)-\nabla_{g}\big(\frac{1}{2}\rho^2\big)\big(\frac{1}{2}|\nabla_{g}u|_{g}^2\big)\\
&=\frac{\partial}{\partial x_{i}}(\rho
   a^{lm}\rho_{x_{m}}u_{x_{l}}a^{ij}u_{x_{j}})-\frac{1}{2}D^2\rho^2\big(\nabla_{g}u,~
\nabla_{g}u\big)-\rho
a^{lm}\rho_{x_{m}}\frac{\partial}{\partial_{x_{l}}}\big(\frac{1}{2}|\nabla_{g}u|_{g}^2\big)\\
&=\frac{\partial}{\partial x_{i}}(\rho
   a^{lm}\rho_{x_{m}}u_{x_{l}}a^{ij}u_{x_{j}})-\frac{1}{2}D^2\rho^2(\nabla_{g}u,~
   \nabla_{g}u)\\
   &-\frac{1}{2}\frac{\partial}{\partial x_{l}}(a^{ij}u_{x_{i}}u_{x_{j}}\rho
   a^{lm}\rho_{x_{m}})+\frac{1}{2}\frac{\partial}{\partial x_{l}}(\rho
   a^{lm}\rho_{x_{m}})|\nabla_{g}u|_{g}^2,
\end{aligned}
 \end{equation}
 and
\begin{eqnarray*}
 \rho a^{lm}\rho_{x_{m}}u_{x_{l}}u^5
   & = & \frac{1}{6}\Big(\rho a^{lm}\rho_{x_{m}}\frac{\partial u^6}{\partial x_{l}}\Big)\\
   & = & \frac{1}{6}\Big[\frac{\partial}{\partial x_{l}}(\rho a^{lm}\rho_{x_{m}}u^6)-u^6\frac{\partial}{\partial x_{l}}
   (\rho a^{lm}\rho_{x_{m}})\Big].
 \end{eqnarray*}
 Finally,
  \begin{eqnarray*}
 0
   & = & u\Big(u_{tt}-\frac{\mathrm{\partial}}{\mathrm{\partial}x_{i}}\big(a^{ij}(x)u_{x_{j}}\big)+u^5\Big)\\
   & = & (uu_{t})_{t}-u_{t}^2-\frac{\mathrm{\partial}}{\mathrm{\partial}x_{i}}\big(ua^{ij}u_{x_{j}}\big)
   +u_{x_{i}}a^{ij}(x)u_{x_{j}}+u^6\\
   & = & \frac{\partial}{\partial t}\big(uu_{t}\big)-\frac{\partial}{\partial
   x_{i}}\big(ua^{ij}u_{x_{j}}\big)+|\nabla_{g}u|_{g}^2+u^6-u_{t}^2.
 \end{eqnarray*}
 Adding, we obtain that
\begin{equation}\label{eq4}
 \begin{aligned}
&\big(tu_{t}+\rho
a^{lm}\rho_{x_{m}}u_{x_{l}}+u\big)\big(u_{tt}-\frac{\mathrm{\partial}}{\mathrm{\partial}x_{i}}(a^{ij}(x)u_{x_{j}})+u^5\big)\\
=&\partial_{t}(tQ+u_{t}u)-\frac{\partial}{\partial x_{i}}(tP)+R=0,
\end{aligned}
 \end{equation}
 where
\begin{equation}
\begin{aligned}
&Q=\frac{1}{2}\Big(u_{t}^2+a^{ij}(x)u_{x_{i}}u_{x_{j}}+\frac{u^6}{3}\Big)+\frac{u_{t}\rho
a^{ij}\rho_{x_{j}}u_{x_{i}}}{t},\nonumber\\
&P=\frac{\rho
a^{ij}\rho_{x_{j}}}{t}\big[\frac{1}{2}\big(u_{t}^2-a^{lm}u_{x_{l}}u_{x_{m}}-\frac{u^6}{3}\big)\big]+
a^{ij}u_{x_{j}}\big(u_{t}+\frac{\rho
a^{lm}\rho_{x_{m}}u_{x_{l}}}{t}+\frac{u}{t}\big),\nonumber\\
&R=\Big(\frac{1}{2}\frac{\partial}{\partial x_{i}}\big(\rho
a^{ij}\rho_{x_{j}}\big)-\frac{3}{2}\Big)u_{t}^2+\frac{1}{2}D^2\rho^2(\nabla_{g}u,~
\nabla_{g}u)\\
&+\Big(\frac{1}{2}-\frac{1}{2}\frac{\partial}{\partial
x_{i}}\big(\rho
a^{ij}\rho_{x_{j}}\big)\Big)|\nabla_{g}u|_{g}^2+\Big(\frac{5}{6}-\frac{1}{6}\frac{\partial}{\partial
x_{i}}\big(\rho a^{ij}\rho_{x_{j}}\big)\Big)u^6.\nonumber
\end{aligned}
\end{equation}
\indent Note that the boundary of the truncated cones ~$K_{S}^T$~ is
 \begin{equation}
 \begin{aligned}
\partial K_{S}^T=\big(([S, ~T]\times \partial\Omega)\cap
K_{S}^T\big)\cup M_{S}^T \cup D(T) \cup D(S),\nonumber\\
\end{aligned}
 \end{equation}
 for~$K_{S}^T,~M_{S}^T,~D(T)$~are as above. Thus if~$\nu_{\partial\Omega}$~denotes the outward unit normal for
 ~$\Omega$,~integrating the identity~$\eqref{eq4}$~over the
truncated
 cones~$K_{S}^T$,~we get
\begin{equation}\label{eq5}
 \begin{aligned}
0=&\int_{D(T)}(TQ+u_{t}u)\mathrm{d}x-\int_{D(S)}(SQ+u_{t}u)\mathrm{d}x+\int_{M_{S}^T}\frac{tQ+uu_{t}-tP\cdot
\nabla\rho}{\sqrt{1+|\nabla\rho|^2}}\mathrm{d}\sigma\\
&+\int_{([S, T]\times \partial\Omega)\cap
K_{S}^T}\nu_{\partial\Omega}\cdot(-tP)\,\mathrm{d}\sigma+\int_{K_{S}^T}R\mathrm{d}t\mathrm{d}x.\\
\end{aligned}
 \end{equation}
First we compute the second to the last term. Note that on $[S,
T]\times
 \partial\Omega$,
 \begin{equation}
 \begin{aligned}
\nabla_{x}u=(\partial_{\nu}u)\nu,~~~u=u_{t}=0.\nonumber\\
\end{aligned}
 \end{equation}
 Thus,
  \begin{equation}
 \begin{aligned}
P=\frac{\rho
a^{lm}\rho_{x_{m}}}{t}(-\frac{1}{2}a^{ij}u_{x_{i}}u_{x_{j}})+
a^{lm}u_{x_{m}}\frac{\rho a^{ij}\rho_{x_{j}}u_{x_{i}}}{t},\nonumber
\end{aligned}
 \end{equation}
  \begin{equation}
 \begin{aligned}
\nabla_{g}u=a^{ij}u_{x_{j}}=(\partial_{\nu}u)a^{ij}\nu_{j},\nonumber
\end{aligned}
 \end{equation}
so
  \begin{equation}
 \begin{aligned}
P&=\frac{\rho
a^{lm}\rho_{x_{m}}}{t}\big(-\frac{1}{2}(\partial_{\nu}u)a^{ij}\nu_{j}(\partial_{\nu}u)\nu_{i}\big)+
(\partial_{\nu}u)a^{lm}\nu_{m}\frac{\rho
a^{ij}\rho_{x_{j}}(\partial_{\nu}u)\nu_{i}}{t}\\
&=\frac{\rho
a^{lm}\rho_{x_{m}}}{t}\big(-\frac{1}{2}(\partial_{\nu}u)^2a^{ij}\nu_{i}\nu_{j}\big)
+a^{lm}\nu_{m}(\partial_{\nu}u)^2\frac{\rho
a^{ij}\rho_{x_{j}}\nu_{i}}{t}.\nonumber
\end{aligned}
 \end{equation}
Finally we get
  \begin{equation}
 \begin{aligned}
-t\nu\cdot P&=\frac{1}{2}\rho
a^{lm}\rho_{x_{m}}\nu_{l}(\partial_{\nu}u)^2a^{ij}\nu_{i}\nu_{j}-a^{lm}\nu_{l}\nu_{m}
(\partial_{\nu}u)^2\rho a^{ij}\rho_{x_{j}}\nu_{i}\\
&=-\frac{1}{2}a^{lm}\nu_{l}\nu_{m} (\partial_{\nu}u)^2\rho
a^{ij}\rho_{x_{j}}\nu_{i}=-\frac{1}{2}a^{lm}\nu_{l}\nu_{m} (\partial_{\nu}u)^2\rho
\nabla_{g}\rho \cdot \nu.\nonumber\\
\end{aligned}
 \end{equation}
 However, for~$x\in \partial\Omega$,~given that~$x_{0}=0\in
 \partial\Omega$,~we have
  \begin{equation}
 \begin{aligned}
\nabla_{g}\rho(x)=\overrightarrow{T}+\mathcal
{O}(x),~~~~\nu(x)=\nu(0)+\mathcal {O}(x),\nonumber\\
\end{aligned}
 \end{equation}
 where~$\overrightarrow{T}$~is a unit vector tangent
 to~$\partial\Omega$~at~$x_{0}=0$.~Consequently, as~$\nu(0)\cdot \overrightarrow{T}=0$,~
  \begin{equation}
 \begin{aligned}
\nabla_{g}\rho(x)\cdot \nu(x)=\mathcal {O}(|x|^{2})=\mathcal {O}(\rho^{2}),
~~~for~~x \in \partial\Omega.\nonumber\\
\end{aligned}
 \end{equation}
 So the second to the last term in~\eqref{eq5}~is bounded(using lemma
 2.9) by
 \begin{equation}
\begin{aligned}
\sup\limits_{x\in K_{S}^{0}}\rho^{2}\times \int_{(-1,~0)\times
\partial\Omega}\big(\frac {\partial u}{\partial
\nu}\big)^{2}\,\mathrm{d}\sigma(x)\,\mathrm{d}t
\leq C|S|^{2}E(u).\nonumber\\
\end{aligned}
\end{equation}
 For the first term
\begin{equation}
\begin{aligned}
\int_{D(T)}uu_{t}\mathrm{d}x&\leq
\Big(\int_{D(T)}u^6\mathrm{d}x\Big)^\frac{1}{6}\Big(\int_{D(T)}u_{t}^2\mathrm{d}x\Big)^\frac{1}{2}
\Big(\int_{D(T)}1\mathrm{d}x\Big)^\frac{1}{3}\\
&\leq
C|T|\big(E(u,D(T))\big)^\frac{1}{6}\big(E(u,D(T))\big)^\frac{1}{2},\nonumber
\end{aligned}
\end{equation}
and
\begin{equation}
\begin{aligned}
\big|\int_{D(T)}TQ\mathrm{d}x\big|&\leq
\int_{D(T)}\big|\frac{T}{2}\big(u_{t}^2+a^{ij}(x)u_{x_{i}}u_{x_{j}}+\frac{u^6}{3}\big)+\frac{Tu_{t}\rho
a^{ij}\rho_{x_{j}}u_{x_{i}}}{t}\big|\mathrm{d}x\\
&\leq C|T|\big(E(u,D(T))\big)+|T|\int_{D(T)}|u_{t}a^{ij}\rho_{x_{j}}u_{x_{i}}|\mathrm{d}x\\
&\leq
C|T|\big(E(u,D(T))\big)+|T|\int_{D(T)}\frac{u_{t}^2}{2}+\frac{a^{lm}\rho_{x_{l}}\rho_{x_{m}}a^{ij}u_{x_{i}}u_{x_{j}}}
{2}\mathrm{d}x\\
&\leq C_{1}|T|\big(E(u,D(T))\big)\\
&\rightarrow 0,\nonumber
\end{aligned}
\end{equation}
as~$T \rightarrow 0$.~\\
\indent So if~$T \rightarrow 0$,~
\begin{equation}
\begin{aligned}
\int_{D(T)}(TQ+u_{t}u)\mathrm{d}x \rightarrow 0.\nonumber
\end{aligned}
\end{equation}
Let~$T\rightarrow 0$~in the identity~\eqref{eq5},~we conclude that
\begin{equation}\label{eq9}
 \begin{aligned}
&-\int_{D(S)}(SQ+u_{t}u)\mathrm{d}x+\int_{M_{S}^0}\frac{tQ+uu_{t}-tP\cdot
\nabla\rho}{\sqrt{1+|\nabla\rho|^2}}\mathrm{d}\sigma\\
&+\int_{([S, 0]\times \partial\Omega)\cap
K_{S}^0}\nu_{\partial\Omega}\cdot(-tP)\,\mathrm{d}\sigma=-\int_{K_{S}^0}R\mathrm{d}t\mathrm{d}x.
\end{aligned}
 \end{equation}
 Let
 \begin{equation}
 \begin{aligned}
&I=-\int_{D(S)}(SQ+u_{t}u)\mathrm{d}x,\\
&II=\int_{M_{S}^0}\frac{tQ+uu_{t}-tP\cdot
\nabla\rho}{\sqrt{1+|\nabla\rho|^2}}\mathrm{d}\sigma.\nonumber
\end{aligned}
 \end{equation}
 On the surface~$\rho=-t$,~we have
\begin{equation}
 \begin{aligned}
&tQ+uu_{t}-tP\cdot \nabla\rho\\
=&\frac{t}{2}\big(u_{t}^2+a^{ij}(x)u_{x_{i}}u_{x_{j}}+\frac{u^6}{3}\big)+u_{t}\rho a^{ij}\rho_{x_{j}}u_{x_{i}}+uu_{t}\\
-&\rho
a^{lm}\rho_{x_{l}}\rho_{x_{m}}\big[\frac{1}{2}\big(u_{t}^2-a^{ij}(x)u_{x_{i}}u_{x_{j}}-\frac{u^6}{3}\big)\big]-
a^{ij}u_{x_{j}}\rho_{x_{i}}(u+tu_{t}+\rho
a^{lm}\rho_{x_{m}}u_{x_{l}})\\
=&-\rho u_{t}^2+2u_{t}\rho a^{ij}\rho_{x_{j}}u_{x_{i}}-\frac{(\rho
 a^{ij}\rho_{x_{j}}u_{x_{i}})^2}{\rho}-\frac{u\rho a^{ij}u_{x_{j}}\rho_{x_{i}}}{\rho}+uu_{t}\\
=&-\rho\big(\frac{\rho
a^{ij}\rho_{x_{j}}u_{x_{i}}}{\rho}-u_{t}\big)^2-u\big(\frac{\rho
a^{ij}\rho_{x_{j}}u_{x_{i}}}{\rho}-u_{t}\big).\nonumber
\end{aligned}
 \end{equation}
 If we parameterize~$M_{S}^0$~by
\begin{equation}
\begin{aligned}
\Omega \ni y \rightarrow \big(-\rho(y), ~y\big),~~~\rho\leq
|S|,\nonumber\\
\end{aligned}
\end{equation}
and let~$\upsilon(y)=u(-\rho(y), y)$,~then
~$\mathrm{d}\sigma=\sqrt{1+|\nabla \rho|^2}\mathrm{d}y$,~and
\begin{equation}
\begin{aligned}
\nabla \upsilon=u_{t}(-\nabla \rho)+\nabla u,\nonumber
\end{aligned}
\end{equation}
furthermore
\begin{equation}
\begin{aligned}
a^{ij}\rho_{y_{j}}\upsilon_{y_{i}}&=-a^{ij}\rho_{y_{j}}\rho_{y_{i}}u_{t}+a^{ij}\rho_{x_{j}}u_{x_{i}}\\
&=-u_{t}+a^{ij}\rho_{x_{j}}u_{x_{i}},\nonumber
\end{aligned}
\end{equation}
so
\begin{equation}
 \begin{aligned}
II&=\int_{M_{S}^0}\frac{tQ+uu_{t}-tP\cdot
\nabla\rho}{\sqrt{1+|\nabla\rho|^2}}\mathrm{d}\sigma\\
&=-\int_{M_{S}^0}\frac{1}{\sqrt{1+|\nabla\rho|^2}}\Big[\rho\big(\frac{\rho
a^{ij}\rho_{x_{j}}u_{x_{i}}}{\rho}-u_{t}\big)^2+u\big(\frac{\rho
a^{ij}\rho_{x_{j}}u_{x_{i}}}{\rho}-u_{t}\big)\Big]\mathrm{d}\sigma\\
&=-\int_{\{y\in \Omega:~\rho \leq |S|\}}\rho
\big(a^{ij}\rho_{y_{j}}\upsilon_{y_{i}}\big)^2+\upsilon\big(a^{ij}\rho_{y_{j}}\upsilon_{y_{i}}\big)\mathrm{d}y\\
&=-\int_{\{y\in \Omega:~\rho \leq |S|\}}\frac{1}{\rho}|\rho
a^{ij}\rho_{y_{j}}\upsilon_{y_{i}}+\upsilon|^2\mathrm{d}y+\int_{\{y\in
\Omega:~\rho \leq |S|\}}\frac{\upsilon^2}{\rho}+\frac{\upsilon\rho
a^{ij}\rho_{y_{j}}\upsilon_{y_{i}}}{\rho}\mathrm{d}y.\nonumber\\
\end{aligned}
 \end{equation}
 For
 \begin{equation}
 \begin{aligned}
&\int_{\{y\in \Omega:~\rho \leq |S|\}}\frac{\upsilon\rho
a^{ij}\rho_{y_{j}}\upsilon_{y_{i}}}{\rho}\mathrm{d}y\\
=&\int_{\{y\in \Omega:~\rho \leq |S|\}}a^{ij}\rho_{y_{j}}\frac{\partial}{\partial_{y_{i}}}(\frac{1}{2}\upsilon^2)\mathrm{d}y\\
=&\int_{\{y\in \Omega:~\rho \leq
|S|\}}\frac{\partial}{\partial_{y_{i}}}\big(a^{ij}\rho_{y_{j}}\frac{1}{2}\upsilon^2\big)-
\frac{1}{2}\upsilon^2\frac{\partial}{\partial_{y_{i}}}\big(a^{ij}\rho_{y_{j}}\big)\mathrm{d}y\\
=&\int_{\rho=|S|}\frac{a^{ij}\rho_{y_{j}}\rho_{y_{i}}}{|\nabla
\rho|}\frac{1}{2}\upsilon^2\mathrm{d}s-\int_{\{y\in \Omega:~\rho
\leq
|S|\}}\frac{1}{2}\upsilon^2\frac{\partial}{\partial_{y_{i}}}\big(\frac{\rho
a^{ij}\rho_{y_{j}}}{\rho}\big)\mathrm{d}y\\
=&\int_{\rho=|S|}\frac{u^2}{2|\nabla \rho|}\mathrm{d}s-\int_{\{y\in
\Omega:~\rho \leq
|S|\}}\frac{\upsilon^2}{2\rho}\big[\frac{\partial}{\partial_{y_{i}}}\big(\rho
a^{ij}\rho_{y_{j}}\big)-1\big]\mathrm{d}y,\nonumber
\end{aligned}
 \end{equation}
 where~$\mathrm{d}s$~is the induced Lebesgue measure on the surface~$\rho=-t$.~And so we have
\begin{equation}\label{eq(10)}
 \begin{aligned}
II&=-\int_{\{y\in \Omega:~\rho \leq |S|\}}\frac{1}{\rho}\big|\rho
a^{ij}\rho_{y_{j}}\upsilon_{y_{i}}+\upsilon\big|^2\mathrm{d}y+\int_{\{y\in
\Omega:~\rho \leq |S|\}}\frac{\upsilon^2}{\rho}+\frac{\upsilon\rho
a^{ij}\rho_{y_{j}}\upsilon_{y_{i}}}{\rho}\mathrm{d}y\\
&=-\int_{\{y\in \Omega:~\rho \leq |S|\}}\frac{1}{\rho}\big|\rho
a^{ij}\rho_{y_{j}}\upsilon_{y_{i}}+\upsilon\big|^2\mathrm{d}y+\int_{\rho
= |S|}\frac{u^2}{2|\nabla \rho|}\mathrm{d}s\\
&+\int_{\{y\in \Omega:~\rho \leq
|S|\}}\frac{\upsilon^2}{\rho}\big[\frac{3}{2}-\frac{1}{2}\frac{\partial}{\partial_{y_{i}}}\big(\rho
a^{ij}\rho_{y_{j}}\big)\big]\mathrm{d}y\\
&=\int_{M_{S}^0}\frac{t\big(a^{ij}\rho_{x_{j}}u_{x_{i}}-u_{t}+\frac{u}{\rho}\big)^2}{\sqrt{1+|\nabla
\rho|^2}}\mathrm{d}\sigma+\int_{\rho = |S|}\frac{u^2}{2|\nabla
\rho|}\mathrm{d}s\\
&+\int_{\{y\in \Omega:~\rho \leq
|S|\}}\frac{\upsilon^2}{\rho}\big[\frac{3}{2}-\frac{1}{2}\frac{\partial}{\partial_{y_{i}}}\big(\rho
a^{ij}\rho_{y_{j}}\big)\big]\mathrm{d}y.
\end{aligned}
 \end{equation}
 In ~$D(S)=\{x \in \Omega:\rho(x) < -S\}$,~$t=S$,~and we have
 \begin{equation}
 \begin{aligned}
SQ+uu_{t}=\frac{S}{2}\big(u_{t}^2+a^{ij}u_{x_{i}}u_{x_{j}}+\frac{u^6}{3}\big)+u_{t}\big(u+\rho
a^{ij}\rho_{x_{j}}u_{x_{i}}\big).\nonumber
\end{aligned}
 \end{equation}
 For the second in the right side,using Cauchy-Schwarz inequality in lemma 3.2 we have
\begin{equation}
 \begin{aligned}
u_{t}(u+\rho a^{ij}\rho_{x_{j}}u_{x_{i}})&\leq
|S|\Big[\frac{u_{t}^2}{2}+\frac{(u+\rho
a^{ij}\rho_{x_{j}}u_{x_{i}})^2}{2S^2}\Big]\\
&\leq |S|\Big[\frac{u_{t}^2}{2}+\frac{(u+\rho
a^{ij}\rho_{x_{j}}u_{x_{i}})^2}{2\rho^2}\Big]\\
&=|S|\frac{u_{t}^2}{2}+\frac{|S|}{2}\Big[\frac{u^2}{\rho^2}+(a^{ij}\rho_{x_{j}}u_{x_{i}})^2+\frac{2u\rho
a^{ij}\rho_{x_{j}}u_{x_{i}}}{\rho^2}\Big]\\
&\leq
|S|\frac{u_{t}^2}{2}+\frac{|S|}{2}\Big[\frac{u^2}{\rho^2}+a^{ij}u_{x_{i}}u_{x_{j}}+\frac{2u\rho
a^{ij}\rho_{x_{j}}u_{x_{i}}}{\rho^2}\Big].\nonumber
\end{aligned}
 \end{equation}
 As ~$S < 0$,~we get
 \begin{equation}\label{eq6}
 \begin{aligned}
SQ+uu_{t}\leq \frac{Su^6}{6}-\frac{Su^2}{2\rho^2}-\frac{Su\rho
a^{ij}\rho_{x_{j}}u_{x_{i}}}{\rho^2},\nonumber
\end{aligned}
 \end{equation}
 so
 \begin{equation}\label{eq7}
 \begin{aligned}
I&=-\int_{D(S)}(SQ+u_{t}u)\mathrm{d}x\\
&\geq
-S\int_{D(S)}\frac{u^6}{6}\mathrm{d}x+S\Big(\frac{1}{2}\int_{D(S)}\frac{u^2}{\rho^2}\mathrm{d}x+\int_{D(S)}
\frac{u\rho a^{ij}\rho_{x_{j}}u_{x_{i}}}{\rho^2}\mathrm{d}x\Big)\\
&=|S|\int_{D(S)}\frac{u^6}{6}\mathrm{d}x+S\Big(\frac{1}{2}\int_{D(S)}\frac{u^2}{\rho^2}\mathrm{d}x+\int_{D(S)}
\frac{u\rho a^{ij}\rho_{x_{j}}u_{x_{i}}}{\rho^2}\mathrm{d}x\Big).
\end{aligned}
 \end{equation}
 Similarly we compute
 \begin{equation}\label{eq8}
 \begin{aligned}
&\int_{D(S)}\frac{u\rho
a^{ij}\rho_{x_{j}}u_{x_{i}}}{\rho^2}\mathrm{d}x\\
=&\int_{D(S)}\frac{a^{ij}\rho_{x_{j}}\frac{\partial}{\partial_{x_{i}}}\Big(\frac{u^2}{2}\Big)}{\rho}\mathrm{d}x\\
=&\int_{D(S)}\frac{\partial}{\partial_{x_{i}}}\Big(\frac{a^{ij}\rho_{x_{j}}u^2}{2\rho}\Big)\mathrm{d}x
-\int_{D(S)}\frac{u^2}{2}\frac{\partial}{\partial_{x_{i}}}\Big(\frac{a^{ij}\rho_{x_{j}}}{\rho}\Big)\mathrm{d}x\\
=&\int_{\rho = |S|}\frac{u^2a^{ij}\rho_{x_{j}}\rho_{x_{i}}}{2\rho
|\nabla
\rho|}\mathrm{d}s-\int_{D(S)}\frac{u^2}{2}\frac{\partial}{\partial_{x_{i}}}
\Big(\frac{\rho a^{ij}\rho_{x_{j}}}{\rho^2}\Big)\mathrm{d}x\\
=&\int_{\rho = |S|}\frac{u^2}{2|S| |\nabla \rho|}\mathrm{d}s-
\int_{D(S)}\frac{u^2}{2}\Big[\frac{1}{\rho^2}\frac{\partial}{\partial_{x_{i}}}\big(\rho
a^{ij}\rho_{x_{j}}\big)-\rho
a^{ij}\rho_{x_{j}}\frac{2\rho_{x_{i}}}{\rho^3}\Big]\\
 =&\int_{\rho =
|S|}\frac{u^2}{2|S| |\nabla
\rho|}\mathrm{d}s-\int_{D(S)}\frac{u^2}{2\rho^2}\Big[\frac{\partial}{\partial_{x_{i}}}\Big(\rho
a^{ij}\rho_{x_{j}}\Big)-2\Big]\mathrm{d}x.
\end{aligned}
 \end{equation}
 Combining~$\eqref{eq7}$~and~$\eqref{eq8}$,~we quickly get
 \begin{equation}\label{eq(11)}
 \begin{aligned}
I\geq |S|\int_{D(S)}\frac{u^6}{6}\mathrm{d}x-\int_{\rho =
|S|}\frac{u^2}{2|\nabla
\rho|}\mathrm{d}s+S\int_{D(S)}\frac{u^2}{2\rho^2}\big[3-\frac{\partial}{\partial_{x_{i}}}(\rho
a^{ij}\rho_{x_{j}})\big]\mathrm{d}x.
\end{aligned}
 \end{equation}

By continuity, the sectional curvature is uniformly bounded near
$x_0$. Then following
from~$\eqref{eq9}$,~$\eqref{eq(10)}$,~$\eqref{eq(11)}$,~and lemma
2.8, using Cauchy-Schwarz inequality in lemma 3.2, we have
 \begin{equation}\label{eq12}
 \begin{aligned}
&|S|\int_{D(S)}\frac{u^6}{6}\mathrm{d}x\leq I+\int_{\rho =
|S|}\frac{u^2}{2|\nabla
\rho|}\mathrm{d}s-S\int_{D(S)}\frac{u^2}{2\rho^2}\Big[3-\frac{\partial}{\partial_{x_{i}}}(\rho
a^{ij}\rho_{x_{j}})\Big]\mathrm{d}x\\
&=-II-\int_{K_{S}^0}R\mathrm{d}t\mathrm{d}x+\int_{([S, ~0]\times
\partial\Omega)\cap
K_{S}^0}\nu_{\partial\Omega}\cdot
(tP)\,\mathrm{d}\sigma\\
&+\int_{\rho = |S|}\frac{u^2}{2|\nabla
\rho|}\mathrm{d}s-S\int_{D(S)}\frac{u^2}{2\rho^2}\Big[3-\frac{\partial}{\partial_{x_{i}}}(\rho
a^{ij}\rho_{x_{j}})\Big]\mathrm{d}x\\
&=\int_{M_{S}^0}\frac{|t|(a^{ij}\rho_{x_{j}}u_{x_{i}}-u_{t}+\frac{u}{\rho})^2}{\sqrt{1+|\nabla
\rho|^2}}\mathrm{d}\sigma-\int_{\{y\in \Omega:~\rho \leq
|S|\}}\frac{\upsilon^2}{\rho}\Big[\frac{3}{2}-\frac{1}{2}\frac{\partial}{\partial_{y_{i}}}(\rho
a^{ij}\rho_{y_{j}})\Big]\mathrm{d}y\\
&-S\int_{D(S)}\frac{u^2}{2\rho^2}\Big[3-\frac{\partial}{\partial_{x_{i}}}(\rho
a^{ij}\rho_{x_{j}})\Big]\mathrm{d}x+\int_{([S, ~0]\times
\partial\Omega)\cap
K_{S}^0}\nu_{\partial\Omega}\cdot
(tP)\,\mathrm{d}\sigma-\int_{K_{S}^0}R\mathrm{d}t\mathrm{d}x\\
&\leq
C|S|\int_{M_{S}^0}(u_{t}-a^{ij}\rho_{x_{j}}u_{x_{i}})^2\mathrm{d}\sigma+C\int_{M_{S}^0}|t|\frac{u^2}{\rho^2}
\mathrm{d}\sigma+C\int_{\{y\in \Omega:~\rho \leq
|S|\}}\frac{\upsilon^2}{\rho}\rho\mathrm{d}y\\
&+C|S|\int_{D(S)}\frac{u^2}{2\rho^2}\rho\mathrm{d}x
+C\int_{K_{S}^0}\rho\big(u_{t}^2+a^{ij}u_{x_{i}}u_{x_{j}}+u^6\big)\mathrm{d}t\mathrm{d}x+C|S|^{2}E(u)-\frac{1}{3}
\int_{K_{S}^0}u^6\mathrm{d}t\mathrm{d}x\\
&\leq
C|S|\int_{M_{S}^0}\big(u_{t}^2-2u_{t}a^{ij}\rho_{x_{j}}u_{x_{i}}+(a^{ij}\rho_{x_{j}}u_{x_{i}})^2\big)\mathrm{d}\sigma
+C\int_{M_{S}^0}\frac{u^2}{|t|}\mathrm{d}\sigma+C\int_{M_{S}^0}\frac{u^2}{\sqrt{1+|\nabla\rho|^2}}\mathrm{d}\sigma\\
&+C|S|\Big(\int_{D(S)}u^6\mathrm{d}x\Big)^\frac{1}{3}
\Big(\int_{D(S)}\frac{1}{\rho^\frac{3}{2}}\mathrm{d}x\Big)^\frac{2}{3}+C\int_{K_{S}^0}
\rho\big(u_{t}^2+a^{ij}u_{x_{i}}u_{x_{j}}+u^6\big)
\mathrm{d}t\mathrm{d}x\\
&+C|S|^{2}E(u)-\frac{1}{3}
\int_{K_{S}^0}u^6\mathrm{d}t\mathrm{d}x\\
&\leq
C_{1}|S|\int_{M_{S}^0}(u_{t}^2+a^{ij}u_{x_{i}}u_{x_{j}})\mathrm{d}\sigma+
C\Big(\int_{M_{S}^0}|t|^{-\frac{3}{2}}\mathrm{d}\sigma\Big)^\frac{2}{3}\Big(\int_{M_{S}^0}u^6\mathrm{d}\sigma\Big)^\frac{1}{3}\\
&+C_{2}\Big(\int_{M_{S}^0}1\mathrm{d}\sigma\Big)^\frac{2}{3}\Big(\int_{M_{S}^0}u^6\mathrm{d}\sigma\Big)^\frac{1}{3}
+C_{3}|S|^2\big(\int_{D(S)}u^6\mathrm{d}x\big)^{\frac{1}{3}}\\
&+C\int_{K_{S}^0}\rho\big(u_{t}^2+a^{ij}u_{x_{i}}u_{x_{j}}+u^6\big)\mathrm{d}t\mathrm{d}x+C|S|^{2}E(u)\\
&\leq C_{4}|S|Flux(u, M_{S}^0)+C_{5}(|S|+|S|^2)\Big(Flux(u,
M_{S}^0)\Big)^\frac{1}{3}+C_{6}|S|^2\big(E(u,~D(S))\big)^{\frac{1}{3}}\\
&+C\int_{K_{S}^0}\rho(u_{t}^2+
a^{ij}u_{x_{i}}u_{x_{j}}+u^6)\mathrm{d}t\mathrm{d}x+C|S|^{2}E(u)\\
&\leq C_{4}|S|Flux(u, M_{S}^0)+C_{5}(|S|+|S|^2)\Big(Flux(u,
M_{S}^0)\Big)^\frac{1}{3}+C_{6}|S|^2E_{0}^{\frac{1}{3}}\\
&+C\int_{K_{S}^0}\rho(u_{t}^2+
a^{ij}u_{x_{i}}u_{x_{j}}+u^6)\mathrm{d}t\mathrm{d}x+C|S|^{2}E(u).\\
\end{aligned}
 \end{equation}
 We put some specific computations in the last part of section
 3, such as the term
  \begin{equation}
 \begin{aligned}
\int_{M_{S}^0}|t|^{-\frac{3}{2}}\mathrm{d}\sigma.\nonumber\\
\end{aligned}
 \end{equation}
 Then combing
with~$\eqref{21777}$~and~$\eqref{eq12}$~we have
 \begin{equation}
 \begin{aligned}
&\int_{D(S)}\frac{u^6}{6}\mathrm{d}x\leq CFlux(u,
M_{S}^0)+C(1+|S|)(Flux(u, M_{S}^0))^\frac{1}{3}\\
&+C_{6}|S|E_{0}^{\frac{1}{3}}+\frac{C\int_{K_{S}^0}\rho(u_{t}^2+
a^{ij}u_{x_{i}}u_{x_{j}}+u^6)\mathrm{d}t\mathrm{d}x}{|S|}+C|S|E(u)\\
&\rightarrow 0~~~~~as~~S \rightarrow 0.\nonumber\\
\end{aligned}
 \end{equation}
 which completes the proof of lemma 2.2.
\section{\textbf{Appendix}}
In this section we give some definition and proofs about Riemannian
Geometry.\\
{\bf Definition 3.1.} Distance function\\
\indent Suppose~$(M,~
g)$~is a Riemannian manifold. For~$x, y\in M$,~we define a
function~$d:M \times M \rightarrow [0, ~\infty)$:~
  \begin{equation}
 \begin{aligned}
d(x,~ y)=\inf \{L(\gamma)|~~\gamma ~~is ~a ~piecewise~ smooth
~curve~ joining ~x~ and ~y\}.
\end{aligned}
 \end{equation}
 If~$M$~is connected, the distance~$d(x, ~y)$~is well defined, since
 there are piecewise smooth curves joining~x~and~$y$.~In this
 case, we can see the function~$d$~satisfies the three properties
 of distance.\\
{\bf Lemma 3.2.} Cauchy-Schwarz inequality\\
If~$A$~is a symmetric, nonnegative~$n\times n$~matrix, then for ~$x,
y\in \mathbb{R}^n$~we have
  \begin{equation}
 \begin{aligned}
\big|\sum_{i,j=1}^{n}a^{ij}x_{i}y_{j}\big|\leq
\big(\sum_{i,j=1}^{n}a^{ij}x_{i}x_{j}\big)^\frac{1}{2}\big(\sum_{i,j=1}^{n}a^{ij}y_{i}y_{j}\big)^\frac{1}{2}.
\end{aligned}
 \end{equation}
{\bf Lemma 3.3.} Suppose~$M$~is a Riemannian manifold, and~$O\in
M$,~let
  \begin{equation}
 \begin{aligned}
\rho:~M\rightarrow [0, ~\infty)~~\rho(x)=d(x, ~O),\nonumber
\end{aligned}
 \end{equation}
 then~$\rho^2\in C^\infty(M)$~in a neighborhood of~$O$,~denoted
 by~$U_{O}$.~And in~$U_{O}$~we have
  \begin{equation}
 \begin{aligned}
|\nabla_{g}\rho|_{g}^2=g^{ij}\rho_{x_{i}}\rho_{x_{j}}=1,~~~~D^2\rho^2
>0,\nonumber
\end{aligned}
 \end{equation}
 where~$(g^{ij})=(g_{ij})^{-1}$.~\\
 {\bf Proof of lemma 2.7.} Fist we compute
 \begin{equation}
 \begin{aligned}
X\big(\frac{1}{2}|\nabla_{g}f|_{g}^2\big)&=\frac{1}{2}X<\nabla_{g}f,~
\nabla_{g}f>_{g}\\
&=<\nabla_{X}\nabla_{g}f,~ \nabla_{g}f>_{g}\\
&=<\nabla_{\nabla_{g}f}X, ~\nabla_{g}f>_{g}+<[X, \nabla_{g}f],~
\nabla_{g}f>_{g}\\
&=<\nabla_{\nabla_{g}f}X,~ \nabla_{g}f>_{g}+[X,~ \nabla_{g}f]f\\
&=<\nabla_{\nabla_{g}f}X,~
\nabla_{g}f>_{g}+X\nabla_{g}f(f)-\nabla_{g}fX(f)\\
 &=<\nabla_{\nabla_{g}f}X, ~\nabla_{g}f>_{g}+X<\nabla_{g}f,~
\nabla_{g}f>_{g}-\nabla_{g}f<X, ~\nabla_{g}f>_{g}\\
&=<\nabla_{\nabla_{g}f}X, ~\nabla_{g}f>_{g}+X<\nabla_{g}f,~
\nabla_{g}f>_{g}\\
&-<\nabla_{\nabla_{g}f}X,~ \nabla_{g}f>_{g}-<X,~
\nabla_{\nabla_{g}f}\nabla_{g}f>_{g}.\nonumber\\
\end{aligned}
 \end{equation}
 So we get
 \begin{equation}
 \begin{aligned}
<X,~
\nabla_{\nabla_{g}f}\nabla_{g}f>_{g}=X\big(\frac{1}{2}|\nabla_{g}f|_{g}^2\big)\nonumber\\
\end{aligned}
 \end{equation}
 Hence one can finish the proof as follows
 \begin{equation}
 \begin{aligned}
<\nabla_{g}f,~
\nabla_{g}\big(X(f)\big)>_{g}&=\nabla_{g}f\big(X(f)\big)\\
&=\nabla_{g}f<X,~\nabla_{g}f>_{g}\\
&= <\nabla_{\nabla_{g}f}X,
\nabla_{g}f>_{g}+<X,~\nabla_{\nabla_{g}f}\nabla_{g}f>_{g}\\
&= <\nabla_{\nabla_{g}f}X,
\nabla_{g}f>_{g}+X\big(\frac{1}{2}|\nabla_{g}f|_{g}^2\big).\nonumber\\
\end{aligned}
 \end{equation}
  {\bf Proof of lemma 2.8.} To do this we need some computation and an additional lemma,
 that is lemma 3.4 as below. Let~$G$~denote
 ~$det(g_{ij})$,~first we compute
  \begin{equation}
 \begin{aligned}
\frac{\partial}{\partial_{x_{i}}}\big(\sqrt{G}
g^{ij}\frac{\partial}{\partial_{x_{j}}}(\frac{1}{2}\rho^2)\big)&=\frac{\partial}{\partial_{x_{i}}}(\sqrt{G}\rho
g^{ij}\rho_{x_{j}})=\frac{\partial \sqrt{G}}{\partial x_{i}}\rho
g^{ij}\rho_{x_{j}}+\sqrt{G}\frac{\partial}{\partial
x_{i}}(\rho g^{ij}\rho_{x_{j}})\\
&=\frac{1}{2\sqrt{G}}\frac{\partial G}{\partial x_{i}}\rho
\nabla_{g}\rho+\sqrt{G}\frac{\partial}{\partial_{x_{i}}}(\rho
g^{ij}\rho_{x_{j}}).\nonumber
\end{aligned}
 \end{equation}
 So
  \begin{equation}\label{eq42}
 \begin{aligned}
\frac{\partial}{\partial_{x_{i}}}(\rho
g^{ij}\rho_{x_{j}})&=\frac{1}{\sqrt{G}}\frac{\partial}{\partial_{x_{i}}}\big(\sqrt{G}
g^{ij}\frac{\partial}{\partial_{x_{j}}}(\frac{1}{2}\rho^2)\big)-\frac{1}{2G}\frac{\partial
G}{\partial g_{lm}}\frac{\partial g_{lm}}{\partial
x_{i}}\rho\nabla_{g}\rho\\
&=\triangle_{g}(\frac{1}{2}\rho^2)-\frac{1}{2}g^{lm}\frac{\partial
g_{lm}}{\partial x_{i}}\rho\nabla_{g}\rho,\\
\end{aligned}
 \end{equation}
 where~$\triangle_{g}$~is the Laplace operator on the Riemannian
 manifold~$(M, ~g)$.~\\
{\bf Lemma 3.4.} Suppose~$M$~is a connected Riemannian manifold,
~$x\in M$,~$\rho$~is the distance function from some point to~$x$.~
If the sectional curvature~$\kappa$~of~$M$~satisfies
  \begin{equation}
 \begin{aligned}
-a^2\leq \kappa\leq a^2,\nonumber
\end{aligned}
 \end{equation}
 where~$a$~is a positive real number. Then in ~$M\setminus\{x~\}$,~ we have
  \begin{equation}\label{eq43}
 \begin{aligned}
1+2a\rho \cot a\rho\leq \triangle_{g}(\frac{1}{2}\rho^2)\leq 1+2a\rho \coth a\rho,\\
\end{aligned}
 \end{equation}
  \begin{equation}\label{eq44}
 \begin{aligned}
a\rho\cot a\rho g\leq D^2(\frac{1}{2}\rho^2)\leq a\rho\coth a\rho g.\\
\end{aligned}
 \end{equation}
It is a classical comparison theorem about the Hessian and Laplace
of the distance function, and one can find the proof in many books
about Riemannian Geometry such as Cheeger and Ebin \cite{Cheeger},
 Greene and Wu \cite{Greene}.\\
 \indent Combining~$\eqref{eq42}$~and~$\eqref{eq43}$~,we quickly get
 the identity~$\eqref{eq217}$~in lemma 2.8. And the identity
 ~$\eqref{eq218}$~can be easily obtained from ~$\eqref{eq44}$.~\\
\indent Now we introduce the geodesic polar coordinates. In this
 coordinate system, the metric can described as follows
 \begin{equation}
 \begin{aligned}
ds^2=d\rho^2+\rho^2g_{11}d\theta^2+2\rho^2g_{12}d\theta
d\varphi+\rho^2g_{22}d\varphi^2,\nonumber
\end{aligned}
 \end{equation}
 and then the estimates of the integrate in the
 geodesic ball or on the mantle can be easily get. For example, we
 estimate~$\int_{D(S)}\frac{1}{\rho^\frac{3}{2}}\mathrm{d}x$~and~$\int_{M_{S}^0}|t|^{-\frac{3}{2}}\mathrm{d}\sigma
 $~in identity~$\eqref{eq12}$:~
  \begin{equation}
 \begin{aligned}
\int_{D(S)}\frac{1}{\rho^\frac{3}{2}}\mathrm{d}x&=\int_{0}^{|S|}\int_{0}^{2\pi}\int_{-\pi}^{\pi}
\frac{\rho^2\sqrt{G}}{\rho^\frac{3}{2}}\,\mathrm{d}\varphi\,\mathrm{d}\theta\,\mathrm{d}\rho\\
&\leq C\int_{0}^S\rho^{\frac{1}{2}}\,\mathrm{d}\rho\\
&\leq C|S|^{\frac{3}{2}},\nonumber\\
\end{aligned}
 \end{equation}
   \begin{equation}
 \begin{aligned}
\int_{M_{S}^0}|t|^{-\frac{3}{2}}\mathrm{d}\sigma&=\int_{0}^{|S|}|t|^{-\frac{3}{2}}\,\mathrm{d}t
\int_{0}^{2\pi}\int_{-\pi}^{\pi}
|t|^2\sqrt{G}\,\mathrm{d}\varphi\,\mathrm{d}\theta\\
&\leq C\int_{0}^{|S|}\sqrt{t}\,\mathrm{d}t\\
&\leq C|S|^{\frac{3}{2}}.\nonumber\\
\end{aligned}
 \end{equation}
\section*{\textbf{Acknowledgement}}
\indent We are very grateful to Professor Yuxin Dong and Professor
Yuanlong Xin for helping us to understand some knowledge of
Riemannian Geometry.\\
\indent The authors are supported by the National Natural Science
Foundation of China under grant 10728101, the 973 Project of the
Ministry of Science and Technology of China, the doctoral program
foundation of the Ministry Education of China, the "111" project and
SGST 09DZ2272900.

\end{document}